\documentclass[12pt]{amsart} 
\usepackage{amsthm, amscd}
\usepackage{amssymb,latexsym} 
\input xy
\xyoption{all}

\theoremstyle{plain}
\newtheorem*{theoremA}{Theorem A}
\newtheorem*{theoremB}{Theorem B}
\newtheorem*{theoremC}{Theorem C}
\newtheorem{theorem}{Theorem}[section]
\newtheorem{proposition}[theorem]{Proposition}

\newtheorem{lemma}[theorem]{Lemma}
 
\theoremstyle{definition}    
\newtheorem{remark}[theorem]{Remark}

\theoremstyle{remark}

\newcommand{\ZZ}{\mathbf{Z}}
\newcommand{\CC}{\mathbf{C}}

\newcommand{\OO}{\mathcal{O}}
\newcommand{\II}{\mathcal{I}}

\newcommand{\EE}{\mathcal{E}}
\newcommand{\FF}{\mathcal{F}}
\newcommand{\LL}{\mathcal{L}}

\renewcommand{\:}{\colon}
\newcommand{\Hom}{\text{Hom}}
\newcommand{\Quot}{\text{\rm Quot}}
\newcommand{\Grass}{\text{Grass}}
\renewcommand{\Im}{\text{Im}}
\newcommand{\Ker}{\text{Ker}}
\newcommand{\IP}{\mathbf{P}}

\begin{document}

\title{\bf Effective very ampleness for generalized theta divisors}
\author{Eduardo Esteves}
\address{\flushleft{\it Eduardo Esteves}\newline 
{\rm Instituto de Matem\'atica Pura e Aplicada, Estrada Dona
Castorina 110, 22460-320 Rio de Janeiro RJ, Brazil; esteves@impa.br}}
\author{Mihnea Popa}
\address{\flushleft{\it Mihnea Popa}\newline 
{\rm Department of Mathematics, Harvard University, One Oxford
Street, Cambridge MA 02138, USA; mpopa@math.harvard.edu}}

\thanks{EE was partially supported by CNPq, Processo 202151/90-5 and 
Processo 300004/95-8, PRONEX, Conv\^enio 41/96/0883/00, and 
FAPERJ, Processo E-26/170.418/2000-APQ1. MP was partially supported by 
NSF Grant DMS-0200150.}

\maketitle
\markboth{EDUARDO ESTEVES AND MIHNEA POPA}
{EFFECTIVE VERY AMPLENESS FOR GENERALIZED THETA DIVISORS}

\section*{\bf Abstract}

\noindent
{\it Given a smooth projective curve $X$, we give effective very ampleness 
bounds for generalized theta divisors on the moduli spaces
$SU_X(r,d)$ and $U_X(r,d)$ of semistable vector bundles of rank $r$ and 
degree $d$ on $X$ with fixed, respectively arbitrary, determinant.}

\section{\bf Introduction}

\noindent
In this paper we address the problem of effective very ampleness for 
linear series on the moduli spaces of vector bundles of arbitrary rank 
on smooth projective curves. We improve and combine techniques used in 
\cite{esteves} and \cite{popa1}, and use the dimension estimate for 
Quot-schemes obtained in \cite{pr}. The problem has been solved before 
only in the case of the determinant line bundle on moduli 
spaces of rank-$2$ vector bundles with fixed determinant, where optimal 
results have been proved in \cite{bv1}, \cite{bv2} and \cite{vgi}. 

Let $X$ be a smooth projective complex curve of genus $g\geq 2$, and fix 
integers $r$ and $d$, with $r>0$. 
Let $U_X(r,d)$ denote the moduli space of equivalence classes 
of semistable vector bundles on $X$ of rank $r$ and degree $d$, and 
$SU_X(r,L)$ the moduli space of semistable rank-$r$ vector bundles
of fixed determinant $L$. The isomorphism class of $SU_X(r,L)$ depends 
only on the degree of $L$, say $d$, so we will often use the notation 
$SU_X(r,d)$ for $SU_X(r,L)$. It is well known (see \cite{dn}, Thm.~B, p.~55) 
that $\text{Pic}(SU_X(r,d))\cong\ZZ$, and the ample generator, henceforth 
denoted by $\LL$, is called the \emph{determinant line bundle}.


\begin{theoremA}\  \newline
For each $m\geq r^2+r$ the linear series $|\LL^m|$ on 
$SU_X(r,d)$ separates points, and is very ample on the 
smooth locus $SU^s_X(r,d)$. In particular, if $r$ and $d$ are coprime, then
$|\LL^m|$ is very ample.
\end{theoremA}

It was already known that $|\LL^m|$ is base-point-free for $m$ in 
the range of the theorem. 
In fact, it is stated in \cite{pr}, Thm.~8.1, p.~661, and proved using the 
method developed in \cite{popa1}, that $|\LL^m|$ is base-point-free for each 
$m\geq [r^2/4]$, at least when $d=0$; see also our Theorem~2.7. 
So the very ampleness bound is roughly four times larger. 
Note that these bounds do not depend on $g$ --- this is important in 
applications, as 
emphasized in \cite{popa1}, especially in Sec.~6.

The lack of knowledge about the tangent space to $SU_X(r,d)$ at a strictly 
semistable point 
prevents us from giving
criteria for the separation of tangent vectors at strictly semistable
points when $r$ and $d$ have a common factor. 

Our approach consists of finding enough generalized theta divisors. These 
divisors arise from vector bundles on $X$ as described below. Set 
$h:=\gcd(r,d)$, $r_1:=r/h$ and $d_1:=d/h$. 
For each vector bundle $F$ on $X$ of rank $mr_1$ and degree 
$m(r_1(g-1)-d_1)$ there is an associated subscheme of $SU_X(r,d)$, 
parameterizing semistable bundles $E$ such that 
$h^0(E\otimes F)\neq 0$. If the subscheme is not the whole 
$SU_X(r,d)$, then it is the support of a divisor, called a 
\emph{generalized theta divisor}, and denoted $\Theta_F$; 
see \cite{dn}, Sec.~0.2, p.~55. The divisor 
$\Theta_F$ belongs to $|\LL^m|$. If $m=1$, the divisor $\Theta_F$ is called 
\emph{basic}.

It was recognized in \cite{faltings} that these divisors can be used to 
produce sections of $\LL^m$ which do not vanish at any given point, 
for $m$ sufficiently large. Later, they 
were also used in \cite{esteves} to separate points and tangent vectors, 
again for $m$ large. However, there is little or no emphasis on effectiveness
in these papers.

Since the generalized theta divisors can be defined by the 
existence of sections or, equivalently, maps between vector bundles,
the most efficient way for producing effective results seems to be a 
counting technique based on dimension estimates for appropriate 
Quot-schemes. The result we need, Lemma 2.1, was proved in \cite{pr}, 
improving on a result of \cite{popa1}.

On $U_X(r,d)$ the situation is slightly different. As above, we can define 
analogous divisors $\Theta_F$, but this time they will only move in the 
same linear series if their determinant, $\det F$, is kept fixed. 
Nevertheless, the same bound is obtained.

\begin{theoremB}\  \newline 
Let $\Theta$ be a basic generalized theta divisor on 
$U_X(r,d)$. For each $m\geq r^2+r$, the linear series 
$|m\Theta|$ on $U_X(r,d)$ separates points, and is 
very ample on the smooth locus $U^s_X(r,d)$. In particular, if $r$ and $d$ 
are coprime, then $|m\Theta|$ is very ample.
\end{theoremB}
Theorem~B is proved using Theorem~A, but is not an immediate corollary. 
Also, methods analogous to those used in the proof of Theorem~A 
would not yield 
\emph{the same} bound on $m$; see the beginning of Section~4. 
We emphasize this since we have encountered a quite widespread, but 
unfounded opinion that such results on $U_X(r,d)$ should follow immediately 
from those on $SU_X(r,d)$ by using the fact that $U_X(r,d)$ admits an 
\'etale cover from $SU_X(r,d)\times J(X)$, where $J(X)$ is the Jacobian of 
$X$. 

To prove Theorem~B we use, as in \cite{popa2}, the Verlinde bundles on 
$J(X)$. This time however we need a more refined cohomological criterion 
for the global generation of arbitrary (as opposed to locally free) 
coherent sheaves on an abelian variety. This result
has been proved recently in \cite{pp} via techniques relying on the
Fourier--Mukai transform.

Finally, combining this technique with the sharp results on rank-2 bundles 
with fixed determinant
mentioned in the first paragraph, we obtain better, and 
presumably optimal, results. 

\begin{theoremC}\  \newline 
Let $\Theta$ be a basic generalized theta divisor 
on $U_X(2,d)$. If $d$ is odd, then $|3\Theta|$ is very ample. 
If $d$ is even, 
and $C$ is not hyperelliptic, 
then $|5\Theta|$ is very ample.
\end{theoremC}

The paper is organized as follows. In Section~2 we prove a very general 
existence result, Lemma~\ref{lemma}, using counting arguments. In Section~3 
we discuss the problem of separating effectively points, 
Theorem~\ref{eff_pts}, 
and tangent vectors, Theorem~\ref{eff_tgv}, and prove 
Theorem~A. Finally, in Section~4 we use Verlinde bundles and 
global generation 
results for coherent sheaves on abelian varieties to show 
Theorems~B~and~C.

\subsection*{\bf Acknowledgments.} The first author 
would like to thank MIT for 
its hospitality while this work was initiated. The second author 
would like to thank 
R. Lazarsfeld and M. Roth for general discussions related to this topic, 
and IMPA for its hospitality and support during part of the work on the 
paper. 

\section{\bf The Existence Lemma}

\noindent
The key computational tool used in the proof of Lemma 2.3, to which this 
section is devoted, is the main 
result of \cite{pr}, stated below. First, a piece of notation: for each 
bundle $E$ on $X$, 
and each $k=0,\dots,r$,
let $d_k(E)$ be the minimum degree of a rank-$k$ 
quotient of $E$.  

\begin{lemma}[\cite{pr}, Thm. 4.1, p. 637]\label{quot}\  \newline
Let $E$ be a vector 
bundle of rank $r$ and degree $d$. Let $\Quot_{k,e}(E)$ be the scheme 
of quotients of $E$ of rank $k$ and degree $e$. Then
	$$\dim\Quot_{k,e}(E)\leq k(r-k)+(e-d_k)r
	\text{ for each $e\ge d_k$, where $d_k:=d_k(E)$.}$$
\end{lemma}

\begin{remark}\  \newline 
At first sight, the appearance of $d_k$ in the bound for 
$\dim\Quot_{k,e}(E)$ might seem undesirable, as it is not clear how 
$d_k$ depends on $E$ and $k$. However, we will see that it is the 
very presence of $d_k$ that allows for the proof of Lemma 2.3.
\end{remark}
 
It is convenient to look at $d_k(E)$ also from a different perspective. 
For each $k$ with $0\leq k\leq r$, put $s_k(E):=rd_{r-k}(E)-(r-k)d$, 
where $r$ and $d$ are the rank and degree of $E$.
(An equivalent and more usual definition is $s_k(E):=kd-rf_k(E)$, 
where $f_k(E)$ is the maximum degree of a rank-$k$ subbundle of $E$; see 
\cite{lange}, p. 450. Also, the usual definition does not contemplate 
the cases $k=0$ and $k=r$.)
Note that $\gcd(r,d)$ divides $s_k(E)$. Also, $E$ is semistable if and only 
if $s_k(E)\geq 0$ for every $k$.


\begin{lemma}\label{lemma}\  \newline 
Let $E$ be a vector bundle of rank $r$ and degree $d$. Fix 
$\epsilon\in\ZZ$. Let $m\in\ZZ$ such that
	$$m\geq\max\big(h,\big[(r+\epsilon)^2/4\big]\big).$$
Let $V$ be the general stable bundle of rank $mr_1$ and degree 
$m(r_1(g-1)-d_1)+\epsilon$. 
If $k$ is the rank of a map $V^*\to E$, then $s_k(E)\le 0$, with equality 
only if $k=0$ or $k<r+\epsilon$.
\end{lemma}

{\noindent
(Note that, if $\epsilon\leq 0$, the lemma says that there are no maps 
$V^*\to E$ of rank $r$.)}

\begin{proof}\  \newline 
We will count moduli. Let 
	$$U:=U_X^s(mr_1,mr_1(g-1)-md_1+\epsilon).$$
Recall that $\dim U=(mr_1)^2(g-1)+1$. For each integers $k$ and $\ell$, 
with $k>0$, put:
	$$U_{k,\ell}:=\{[V]\in U\,|\,\text{there is a map $\psi:V^*\to E$ 
	with image of rank $k$ and degree $\ell$}\}.$$
The lemma will follow once we prove that 
	$$\dim U_{k,\ell}\leq (mr_1)^2(g-1)
	\text{ if $s_k(E)>0$, or if $s_k(E)=0$ and 
	$k\geq r+\epsilon$}.$$

So fix $k$ and $\ell$, with $k>0$. Let $V$ run among the stable bundles of 
rank $mr_1$ and degree $mr_1(g-1)-md_1+\epsilon$, and $\psi:V^*\to E$ among 
the maps whose image has rank $k$ and degree $\ell$. 
Let $F:=\Im(\psi)$ and $G:=\Ker(\psi)$. 
In other terms, 
we have the following diagram of maps, where the horizontal sequence is 
exact, the vertical map is an inclusion, and the maps in the triangle 
commute:
	$$\xymatrix{
	0 \ar[r] & G \ar[r] & V^* \ar[r] \ar[dr]_{\psi} & F \ar[r] 
	\ar@{^{(}->}[d]  & 0\\
	&  &  & E}$$

The bundles $F$, or rather the quotients $E/F$, are parameterized by the 
scheme $Q$ of quotients of $E$ of degree $d-\ell$ and rank $r-k$. 
By Lemma 2.1, letting $d_{r-k}:=d_{r-k}(E)$,
$$
\dim Q\leq k(r-k)+(d-\ell-d_{r-k})r.
$$

The bundles $G$ are subbundles of fixed rank, equal to $mr_1-k$, 
and fixed degree of the stable bundles $V^*$. So they are parameterized by 
a scheme $T$ --- a subset of a relative Quot-scheme over an \'etale cover of 
$U$. For counting purposes, as is probably well known, and observed 
in \cite{bgn}, Rmk. 4.2(i), p. 656, we may assume that $\dim T$ is 
at most the dimension of the moduli space of stable bundles of same rank, 
$mr_1-k$. So
	$$\dim T\leq (mr_1-k)^2(g-1)+1.$$

The bundles $V^*$ are extensions of the bundles $F$ by the bundles $G$. 
These extensions are parameterized by the projectivization of  
$H^1(F^*\otimes G)$. Not all extensions of $F$ by $G$ are stable. In 
fact, none will be unless $\Hom(F,G)=0$, as a stable bundle admits only 
the trivial endomorphisms. If $\Hom(F,G)=0$ then, by the Riemann--Roch 
Theorem, 
	$$h^1(F^*\otimes G)=-\chi(F^*\otimes G)=
	-(mr_1-k)\chi(F^*)-k\deg G=p+1,$$
where
	$$p:=(mr_1-k)(\ell+k(g-1))+k(mr_1(g-1)-md_1+\epsilon+\ell)-1.$$
Let $A$ be the open subset of $Q\times T$ parameterizing the pairs 
$(F,G)$ satisfying $h^0(F^*\otimes G)=0$. Then the 
bundles $V^*$ are parameterized by a projective bundle $P$ over $A$ of 
relative dimension $p$.

Clearly, $\dim U_{k,\ell}\leq\dim P$. Now, 
	\begin{align*}
	\dim P=&\dim(P/A)+\dim Q+\dim T\\
	\leq&(mr_1-k)(\ell+k(g-1))+k(mr_1(g-1)-md_1+\epsilon+\ell)-1\\
	&+k(r-k)+(d-\ell-d_{r-k})r+(mr_1-k)^2(g-1)+1\\
	=&(mr_1)^2(g-1)+k(r+\epsilon-k)+mr_1\ell-kmd_1+(d-\ell-d_{r-k})r.
	\end{align*}
Using $s_k(E)=kd-r(d-d_{r-k})$, we have
	$$mr_1\ell-kmd_1+(d-\ell-d_{r-k})r=
	-r_1(d-\ell-d_{r-k})(m-h)-ms_k(E)/h.$$
This gives
	\begin{equation}\label{dimU}\begin{aligned}
	\dim U_{k,\ell}\leq&(mr_1)^2(g-1)+k(r+\epsilon-k)\\
	&-r_1(d-\ell-d_{r-k})(m-h)-ms_k(E)/h.
	\end{aligned}\end{equation}

Now, $d-\ell\geq d_{r-k}$, since $d-\ell$ is the degree of the quotient 
$E/F$, which has rank $r-k$. Since $m\geq h$, we get
	$$\dim U_{k,\ell}\leq(mr_1)^2(g-1)+k(r+\epsilon-k)-ms_k(E)/h.$$
If $s_k(E)>0$ then $s_k(E)\geq h$, and hence 
$k(r+\epsilon-k)-ms_k(E)/h\leq 0$ from the hypothesis on $m$. 
The same holds if $s_k(E)=0$ and $k\geq r+\epsilon$.
In either case, $\dim U_{k,\ell}\leq(mr_1)^2(g-1)$, as asserted.
\end{proof}

\begin{remark}\  \newline 
As it can be seen from the last three lines in the above proof, it is enough 
to assume that $m\geq h$ and $m\geq hi(r+\epsilon-i)/s_i(E)$ for every $i$ 
such that $s_i(E)>0$ for the conclusion of Lemma~\ref{lemma} to hold. 
In addition, if $E$ is semistable, and if 
$m\geq h+i(r+\epsilon-i)/r_1$ for every $i$ such that $s_i(E)=0$, then 
the image of each map $V^*\to E$ is semistable of slope $d/r$. Indeed, let 
$\phi\:V^*\to E$ be a map, say of rank $k$. Then $s_k(E)\geq 0$ 
because $E$ is semistable, and $s_k(E)\leq 0$ by Lemma~\ref{lemma}. 
So $s_k(E)=0$. Assume, by contradiction, that $F:=\Im(\phi)$ is not 
semistable of slope $d/r$. Since $s_k(E)=0$ and $E$ is semistable, the 
degree $\ell$ of $F$ satisfies $d-\ell>d_{r-k}$. 
Now, $m\geq h+k(r+\epsilon-k)/r_1$ by 
hypothesis, since $s_k(E)=0$. So Formula~(\ref{dimU}) yields 
$\dim U_{k,\ell}\leq(mr_1)^2(g-1)$. Since $[V]\in U_{k,\ell}$, the bundle 
$V$ is not general, a contradiction. 
\end{remark}

\begin{remark}\  \newline
The upper bound for $\dim U_{k,\ell}$ in 
Formula~(\ref{dimU}) may not be sharp in particular cases, as we are 
using a bound on the dimension of the space of quotients, rather than on the 
number of moduli. However, it does not seem possible to improve the bound 
in a uniform way; see \cite{pr}, \S 8.
\end{remark}

\begin{proposition}\label{general_vanishing}\  \newline
Let $E$ be a stable bundle of rank $r$ and degree $d$. 
If $m\geq\max([r^2/4],h)$, the general stable bundle $V$ of rank 
$mr_1$ and degree $m(r_1(g-1)-d_1)$ satisfies $H^0(V\otimes E)=0$. 
\end{proposition}

\begin{proof}\  \newline
Let $k$ be the rank of a map $V^*\to E$. Apply 
Lemma~\ref{lemma} with $\epsilon:=0$. Then $s_k(E)\leq 0$, 
with equality only if $k<r$. However, 
since $E$ is stable, $s_i(E)>0$ for each $i$ with $0<i<r$, and $s_r(E)=0$.
So $k=0$, whence $H^0(V\otimes E)=0$.
\end{proof}

Note that $[r^2/4]\geq r\geq h$ if $r\geq 4$. 


\begin{theorem}[cf. \cite{pr}, Thm. 8.1, p. 661]\  \newline
The linear series $|\LL^m|$ on $SU_X(r,d)$ has no base points for 
each $m\geq\max(h,[r^2/4])$.
\end{theorem}

\begin{proof}\  \newline
Let $E$ be a direct sum of stable bundles $E_i$ 
corresponding to a point $[E]$ of $SU_X(r,d)$. 
Let $k_i$ and $e_i$ be the rank and degree of $E_i$, and set  
$t_i:=\gcd(k_i,e_i)$. Since $e_i/k_i=d/r$, we have $r_1=k_i/t_i$ and 
$d_1=e_i/t_i$. Also, $t_i\leq h$. 

Apply Proposition~\ref{general_vanishing} to each $E_i$. Since 
$m\geq\max(t_i,[k_i^2/4])$, the general stable bundle $V$ of rank 
$mr_1$ and degree $m(r_1(g-1)-d_1)$ satisfies $H^0(V\otimes E_i)=0$. 
So $H^0(V\otimes E)=0$, and thus 
$\Theta_V$ is a divisor of $|\LL^m|$ that does not contain $[E]$ 
in its support. 
\end{proof}

\section{\bf Very ampleness on $SU_X(r,d)$}

\noindent
The existence of generalized theta divisors gives a very simple 
sufficient criterion for a pluritheta linear series to separate points.

\begin{lemma}\label{sep_points}\  \newline
Let $E_1$ and $E_2$ be semistable bundles representing distinct
points of $SU_X(r,d)$. If there exists a bundle $F$ of rank $mr_1$ and 
degree $m(r_1(g-1)-d_1)$ such that $h^0(E_1\otimes F)\neq 0$ and 
$h^0(E_2\otimes F)=0$, then the linear series $|\LL^m|$
separates the points corresponding to $E_1$ and $E_2$. 
\end{lemma}

\begin{proof}\  \newline
The statement simply says that $\Theta_F$ is a divisor 
in $|\LL^m|$ passing through the point corresponding to $E_1$, but not
through that corresponding to $E_2$.   
\end{proof} 

For our purposes, an \emph{elementary transformation} of a vector bundle $E$ 
at a point $P$ of $X$ is the vector bundle $E'$ sitting in an 
exact sequence,
	$$0\longrightarrow E'\longrightarrow E \longrightarrow \CC_P
	\longrightarrow 0,$$
where $\CC_P$ is the skyscraper sheaf at $P$. The map $E\to\CC_P$ 
corresponds to a linear functional on $E(P)$, and $E'$ depends 
only on the hyperplane of zeros of this functional. We will use 
elementary transformations to distinguish between maps from a 
given bundle to $E$.

\begin{theorem}\label{eff_pts}\  \newline
The linear series $|\LL^m|$ on $SU_X(r,d)$ separates points 
for each $m\geq r^2+r$.
\end{theorem}

\begin{proof}\  \newline
Let $G_1$ and $G_2$ be two semistable bundles corresponding to distinct 
points of $SU_X(r,d)$. We may assume $G_1$ and $G_2$ are direct sums of 
stable bundles. If they have a common stable summand $K$, say 
$G_1=H_1\oplus K$ and $G_2=H_2\oplus K$, then our problem is reduced to 
that of separating the points corresponding to $H_1$ and $H_2$ in a moduli 
space of bundles of smaller rank, as done in the proof of 
\cite{esteves}, Thm.~16, p.~590. So we may assume $G_1$ and $G_2$ have no 
common stable summand. In particular, there is a stable summand $E_0$ of 
$G_1$ distinct from the stable summands $E_1,\dots,E_n$ of $G_2$. 
To apply Lemma~\ref{sep_points} we need only to find a bundle 
$F$ of rank $mr_1$ and degree $m(r_1(g-1)-d_1)$ such that 
$h^0(E_0\otimes F)\neq 0$ but $h^0 (E_i\otimes F)=0$ for each
$i=1,\dots,n$. 

For each $i$, let $k_i$ denote the rank of $E_i$. 
We first claim that, for $m$ as in the statement of the theorem, 
the general stable bundle $V$ of rank $mr_1$ and degree 
$m(r_1(g-1)-d_1)+1$ satisfies the following two properties:
\begin{enumerate}
\item For each $i=0,\dots,n$, we have $h^1(E_i\otimes V)=0$ and 
$h^0(E_i\otimes V)=k_i$. 
\item For each $i=1,\dots,n$, and each nonzero maps $\phi\:E_0^*\to V$ and 
$\psi\:E_i^*\to V$, the sum $(\phi,\psi)\:E_0^*\oplus E_i^*\to V$ is 
injective with saturated image.
\end{enumerate} 
Indeed, let us check property 1 first. The bundle $W:=V^*\otimes\omega_C$ is 
the general bundle of rank $mr_1$ and degree $m(r_1(g-1)+d_1)-1$. 
Let $k$ be the rank of a map $W^*\to E^*_i$. 
Apply Lemma~\ref{lemma} 
with $\epsilon:=-1$. Then $s_k(E_i)\leq 0$, with equality only if $k<k_i -1$. 
Since $E_i$ is stable, $k=0$. Thus $h^0(W\otimes E_i^*)=0$. By 
Serre Duality, $h^1(V\otimes E_i)=0$. Riemann--Roch yields now 
$h^0(E_i\otimes V)=k_i$.
 
As for property 2, 
apply Lemma~\ref{lemma} with $\epsilon:=1$ to $E_0$, $E_i$ 
and $E_0\oplus E_i$. Then $\phi^*$ has rank $k$ satisfying 
$s_k(E_0)\leq 0$. Since $\phi\neq 0$ and $E_0$ is stable, $\phi^*$ is 
surjective. Likewise, $\psi^*$ is surjective. The same analysis shows that 
the rank of $(\phi,\psi)^*$ is $k_0$, $k_i$,  or $k_0 + k_i$. 
Now, the image of $(\phi,\psi)^*$ surjects onto $E_0$ and $E_i$, because 
$\phi^*$ and $\psi^*$ are surjective. 
So, if the rank of $(\phi,\psi)^*$ were not $k_0 + k_i$, 
we would obtain an isomorphism between $E_0$ and $E_i$. Hence
$(\phi,\psi)^*$ is surjective, and thus
$(\phi,\psi)$ is injective with saturated image. 

Fix a nonzero map $\phi:E_0^*\to V$ and set $Q:=\text{Coker}(\phi)$. By 
property 2 above, $Q$ is a vector bundle. In addition, for each 
nonzero map $\psi:E_i^*\to V$ the composition 
$\widetilde{\psi}:E_i^*\to V\to Q$ is injective with saturated image. 
Fix any point $P$ of $X$. Then the image of $\widetilde{\psi}(P)$ is a 
vector subspace $L_{\psi,i}$ of dimension $k_i$ of $Q(P)$. Letting $\psi$ 
vary, the $L_{\psi,i}$ move in a subvariety $M_i$ of the Grassmannian 
$\Grass(k_i,Q(P))$. Note that $\dim M_i\leq k_i-1$ because 
$h^0(E_i\otimes V)=k_i$ from property 1 above. Let $H$ be a 
hyperplane of $Q(P)$, and let $Z_i$ be the closed subset of 
$\Grass(k_i,Q(P))$ parameterizing subspaces contained in $H$.
The codimension of $Z_i$ is $k_i$. Choosing $H$ general enough, by 
\cite{kl}, Cor. 4, (i), p. 291, the intersection $M_i\cap Z_i$ is empty 
for each $i=1,\dots,n$.

Let $H'$ be the inverse image of $H$ in $V(P)$, and $F$ the elementary 
transformation of $V$ at $P$ corresponding to the hyperplane $H'$. 
Since $H'$ contains the image of $\phi(P)$, the map $\phi$ factors through 
$F$, whence $h^0(E_0\otimes F)\neq 0$. On the other hand, $H'$ does not 
contain the image of $\psi(P)$ for any nonzero map $\psi:E_i^*\to V$. So 
$h^0(E_i\otimes F)=0$ for each $i=1,\dots,n$. Now, $F$ has rank $mr_1$ and 
degree $m(r_1(g-1)-d_1)$, so, by Lemma \ref{sep_points}, 
$|\LL^m|$ separates the points corresponding to $G_1$ and $G_2$. 
\end{proof}

Let $S$ be the spectrum of the algebra of dual numbers 
$A:=\CC[\epsilon]/(\epsilon ^2)$. If $E$ is a stable bundle, 
representing a point $[E]$ of 
$SU_X(r,d)$, the tangent space to $SU_X(r,d)$ at $[E]$ 
can be identified with the space of deformations of $E$ over $S$ 
with constant determinant.

\begin{lemma}\label{sep_tgv}\  \newline
Let $E$ be a stable bundle representing a point $[E]$ of $SU_X(r,d)$. 
The linear series $|\LL^m|$ separates tangent vectors at $[E]$ if, for each 
nontrivial deformation $\EE$ of $E$ over $S$, there exists a 
bundle $F$ of rank $mr_1$ and degree $m(r_1(g-1)-d_1)$ such that 
$h^0(E\otimes F)=1$ and the unique (modulo $\CC$) 
section of $E\otimes F$ 
does not extend over $S$ to a section of $\EE\otimes F$. 
\end{lemma}

\begin{proof}\  \newline
 Let $G:=\omega_C\otimes F^*$. Let 
$\iota\:SU_X(r,d)\to SU_X(r,-d)$ be the involution taking a vector bundle 
to its dual. By Serre Duality, $\iota^*\Theta_G=\Theta_F$. Let 
$f\: S\to SU_X(r,d)$ be the map induced by $\EE$. Then $\iota f$ is induced 
by $\EE^*$. 

By hypothesis, $h^0(E\otimes F)=1$; so $h^1(E^*\otimes G)=1$ by Serre 
Duality. 
By Riemann--Roch, also $h^0(E^*\otimes G)=1$. From the hypothesis, 
the unique (modulo $\CC$) nonzero map $E^*\to G^*\otimes\omega_C$ does 
not extend over $S$ to a map $\EE^*\to G^*\otimes\omega_C\otimes\OO_S$. 
By \cite{esteves}, Lemma 11, p. 582, $(\iota f)^*\Theta_G$ is reduced. 
Then $f^*\Theta_F$ is reduced as well. So $\Theta_F$ is a divisor in 
$|\LL^m|$ containing $[E]$, whose tangent space at $[E]$ does not contain 
the tangent vector corresponding to $\EE$.
\end{proof}

\begin{theorem} \label{eff_tgv}\  \newline
The linear series $|\LL^m|$ on $SU_X(r,d)$ separates 
tangent vectors at stable points for each $m\geq r^2+r$.
\end{theorem}

\begin{proof}\  \newline 
Let $E$ be a stable bundle of rank $r$ and degree $d$. 
Let $\EE$ be a nontrivial deformation of $E$ over $S$. We need only 
to find a bundle $F$ satisfying the conditions set forth in 
Lemma~\ref{sep_tgv}. 

Let $V$ be the general stable bundle of degree 
$m(r_1(g-1)-d_1)+1$ and rank $mr_1$. Let 
$\lambda,\mu\in\Hom(V^*,E)-\{0\}$. 
By Lemma~\ref{lemma}, since $E$ is stable, both $\lambda$ and $\mu$ 
are surjective. Consider the sum $(\lambda,\mu)\:V^*\to E\oplus E$. Again 
by Lemma~\ref{lemma}, and the stability of $E$, 
either $(\lambda,\mu)$ is surjective or it has 
rank $r$. In the latter case, the image of $(\lambda,\mu)$ 
projects isomorphically onto 
both factors, yielding an automorphism $\tau$ of $E$ such that 
$\lambda=\tau\mu$. By stability, $\text{Aut}(E)=\CC^*$. Hence, either 
$(\lambda,\mu)$ is surjective or $\lambda$ is a multiple of $\mu$.

As in the proof of Theorem~\ref{eff_pts}, we also have 
$h^1(E\otimes V)=0$ and $h^0(E\otimes V)=r$.

Fix now a nonzero map $\phi\:E^*\to V$. Since $\phi^*$ is surjective, 
$\phi$ is injective with saturated image. Let $\pi\:V\to W$ be the 
quotient map to its cokernel. Since $h^1(E\otimes V)=0$, the map $\phi$ 
extends over $S$ to an injection $\mu\:\EE^*\to V\otimes\OO_S$ with flat 
cokernel. The injection $\mu$ corresponds to a tangent vector at $[\pi]$
of the Quot-scheme of $V$, whence to a map $\nu\:E^*\to W$. 
All of the extensions of $\phi$ over $S$ correspond to maps of the form 
$\nu+\pi\psi$ for $\psi\in\Hom(E^*,V)$. All of these maps are nonzero 
because $\EE$ is nontrivial.

For each nonzero map $\psi\: E^*\to V$ consider the dual map 
$(\phi,\psi)^*\: V^*\to E\oplus E$. If 
$\psi$ is not a multiple of $\phi$ then $(\phi,\psi)^*$ is surjective, and 
thus $(\phi,\psi)\: E^*\oplus E^*\to V$ is injective with saturated image. 
So either $\pi\psi\:E^*\to W$ is zero, when $\psi$ is a multiple of 
$\phi$, or injective with saturated image. In particular, the vector 
subspace $N\subseteq\Hom(E^*,W)$ generated by $\pi\psi$ for 
$\psi\in\Hom(E^*,V)$ has dimension $r-1$. Note also that $\nu\not\in N$.

Let $\lambda\:N\otimes\OO_X\to Hom(E^*,W)$ be the natural map. 
Since each nonzero map in $N$ is injective with saturated image, 
$\lambda$ is also injective with saturated image. 
Thus, if $\nu(P)$ belonged to $\Im(\lambda(P))$ for each $P\in X$, 
the section of $Hom(E^*,W)$ corresponding to $\nu$ would factor 
through $\lambda$. We would get $\nu\in N$, a contradiction.

So pick a point $P$ of $X$ such that $\nu(P)\not\in\Im(\lambda(P))$. 
Let $A:=\Im(\lambda(P))$ and $B\subseteq\Hom(E^*(P),W(P))$ the subspace 
generated by $A$ and $\nu(P)$. Consider the associated projective spaces 
$\IP(A)$ and $\IP(B)$ of one-dimensional subspaces of $A$ and $B$. 
Then $\IP(B)$ has dimension $r-1$ and $\IP(A)$ is a hyperplane in
$\IP(B)$. For each $Q\in\IP(B)$ denote by $L_Q$ the image of the 
unique (modulo $\CC^*$) map $E^*(P)\to W(P)$ corresponding to $Q$. 
Let $U\subseteq\IP(B)$ be the open subset of points $Q$ such that $L_Q$
has 
maximal dimension $r$. Then $U\supseteq\IP(A)$, and hence the complement
of 
$U$ in $\IP(B)$ is a finite set $T$. For each 
$Q\in T$ the image $L_Q$ does not have maximal dimension, but is nonzero 
because $\nu(P)\not\in A$.

Let $H$ be a general hyperplane of $W(P)$. So $H$ 
does not contain any of the finitely many subspaces 
$L_Q$ for $Q\in T$. Also, let $S$ be the closed subset 
of $\text{Grass}(r,W(P))$ parameterizing subspaces contained in $H$. 
The codimension of $S$ is $r$. So, by \cite{kl}, Cor. 4, (i), p.~291, 
since $H$ is general, 
the intersection of $S$ and 
$\{[L_Q]\,|\,Q\in U\}$ is empty. The upshot is that $H$ does not contain
$L_Q$ 
for any $Q\in\IP(B)$. 

Let $G$ be the elementary transformation of $W$ at $P$ 
corresponding to the hyperplane $H$,
and $F$ its inverse image under 
$\pi\:V\to W$. Then $F$ has rank $mr_1$ and degre $m(r_1(g-1)-d_1)$. 
In addition, $\phi$ factors through an injection 
$\phi'\:E^*\to F$ with quotient map $\pi'\:F\to G$ induced by $\pi$. 
However, those maps $\psi\:E^*\to V$ that are not multiples of $\phi$ 
do not factor 
through $F$, since $H$ does not contain the image of $\pi\psi(P)$ for
any such $\psi$. Hence $h^0(E\otimes F)=1$, and every map $E^*\to F$ 
is a multiple of $\phi'$.

Finally, $\phi'$ does not extend over $S$ to a map 
$\EE^*\to F\otimes\OO_S$. In fact, if $\phi'$ had such an extension, 
coupling it with the inclusion $i\:F\to V$, we would get an extension 
$\EE^*\to V\otimes\OO_S$ of $\phi$ over $S$. So there would be a map 
$\nu'\:E^*\to G$, corresponding to the extension of $\phi'$, whose 
composition with the inclusion $G\subseteq W$ is of the form 
$\nu+\pi\psi$ for a certain $\psi\:E^*\to V$. 
Then $\nu(P)+\pi\psi(P)$ would factor through $H$, reaching a contradiction 
with our choice of $H$. 
\end{proof}

\begin{proof} (of Theorem A)\  \newline 
By Theorem 3.2, the complete linear series $|\LL^m|$ 
corresponds to an injective map $\phi\:SU_X(r,d)\to\text{\bf P}^N$. Since 
$\phi$ is injective and proper, $\phi(SU_X^s(r,d))$ is open in 
$\phi(SU_X(r,d))$ and $SU_X^s(r,d))=\phi^{-1}(\phi(SU^s_X(r,d)))$. 
So the restriction $\psi$ of $\phi$ to $SU_X^s(r,d)$ is injective and 
proper over its image. By Theorem 3.4, $\psi$ is an embedding.
\end{proof}

\section{\bf Very ampleness on $U_X(r,d)$}

\noindent\setcounter{equation}{0}
Let $U_X(r,d)$ be the moduli space 
of semistable bundles of rank $r$ and degree $d$ on $X$. 
Let $m$ be a positive integer. 
As before, to a vector bundle $G$ of rank $mr_1$ 
and degree $m(r_1(g-1)-d_1)$ we can associate a subscheme 
of $U_X(r,d)$ parameterizing semistable bundles $E$ such that 
$h^0(E\otimes G)\neq 0$. If the subscheme is not the whole $U_X(r,d)$, 
then it is the support of a (generalized theta) divisor, denoted 
$\Theta_G$; see \cite{dn}, \S 0.2, p.~55.

We search for effective results --- bounds on $m$ --- towards the very 
ampleness of the linear series $|m\Theta_G|$. 
However, in contrast with the case of $SU_X(r,d)$, 
as $G$ moves among vector bundles of equal rank and degree, 
the divisors $\Theta_G$ will only move in the same linear
series if $\det G$ is fixed; see the formula at the bottom of p.~57 in 
\cite{dn}. As a consequence, the methods presented in Section 3 
do not apply immediately. It is still possible to make them work, by 
considering theta divisors of the form $\Theta_{G_1\oplus G_2}$, where 
each $G_i$ is a bundle of rank $m_ir_1$ and degree 
$m_i(r_1(g-1)-d_1)$, 
with $m_1+m_2=m$. The bundle $G_1$ is chosen so that $\Theta_{G_1}$ has the 
desired separation properties, and $G_2$ is added only to 
have the determinant of $G_1\oplus G_2$ equal to $(\det F)^m$, and thus 
$\Theta_{G_1\oplus G_2}\in|m\Theta_F|$; cf. the 
end of the proof of Lemma 10 on p.~580 in \cite{esteves}. However, this 
procedure would lead to bounds on $m$ worse than those found for 
$SU_X(r,d)$. 

In this section we show that in fact one can obtain for $U_X(r,d)$ 
the same bounds 
as for $SU_X(r,d)$. Our method relies on the use of the Verlinde bundles 
defined in \cite{popa2}. We begin by recalling the notation. 

Let $L$ be a line bundle of degree $d$ on $X$ and $\pi_L$ the composition:
	$$\pi_L\:U_X(r,d)\overset{{\rm det}}{\longrightarrow} {\rm Pic}^d(X)
	\overset{\otimes L^{-1}}{\longrightarrow} J(X).$$
Set $U:=U_X(r,d)$ and $J:=J(X)$. Let
$F$ be a vector bundle on $X$ of rank $mr_1$ and degree 
$m(r_1(g-1)-d_1)$ for which there is a generalized theta divisor 
$\Theta_F$. Put
	$$V_m:={\pi_L}_{*}\OO_U(m\Theta_F).$$
The fibers of $\pi_L$ are the $SU_X(r,A)$, for $A$ running over all line 
bundles of degree $d$ on $X$, and $\OO_U(\Theta_F)_{|SU_X(r,A)}=\LL_A$, the 
determinant line bundle on $SU_X(r,A)$. Since there is an isomorphism 
$SU_X(r,A)\cong SU_X(r,B)$ for each line bundles $A$ and $B$ of degree $d$, 
and the isomorphism carries $\LL_A$ to $\LL_B$, the 
\emph{Verlinde number}, $v_m:=h^0(SU_X(r,A),\LL^m_A)$, does not depend 
on $A$. Hence $V_m$ is a bundle of rank $v_m$, called the 
\emph{Verlinde bundle}.

Consider the diagram of maps,
	$$\xymatrix{
	SU_{X}(r,L)\times J(X) \ar[d]_{p_{2}} \ar[r]^{\hspace{7mm} \tau} & 
	U_{X}(r,d) \ar[d]^{{\pi}_{L}} \\
	J(X) \ar[r]^{r_{J}} & J(X) }$$
where $\tau$ is given by tensor product, $p_2$ is the second projection, 
and $r_J$ is given by multiplication by $r$ on the Jacobian. The diagram is 
clearly commutative. Both $\tau$ and $r_J$ are \'etale coverings 
of degree $r^{2g}$; see \cite{tt}, Prop. 8, p. 338 for $\tau$. Hence the 
diagram is a fiber diagram.

By \cite{dt}, Cor. 6, p. 350,
	\begin{equation}\label{pullback}
	\tau^*\OO_U(\Theta_F)\cong\LL\boxtimes\OO_J(rr_1\Theta_N),
	\end{equation}
where $N$ is any line bundle of degree $g-1$ on $X$ such that 
$N^r\cong L\otimes(\det F)^h$.
Then the base-change formula for the diagram above yields
	\begin{equation}\label{mult_r}
	r_{J}^{*}V_m\cong\bigoplus_{i=1}^{v_m}\OO_{J}(mrr_{1}\Theta_N).
	\end{equation}

Let $\Sigma$ be a finite scheme over $\CC$. A coherent sheaf $\FF$ on a 
scheme $Z$ is called \emph{$\Sigma$-spanned} (resp. \emph{$\Sigma$-spanned 
on an open subset $Y\subseteq Z$}) if for each map $f\:\Sigma\to Z$ 
(resp. $f\:\Sigma\to Y$), each section of $\FF|_{f(\Sigma)}$ lifts to 
a global section of $\FF$. The sheaf $\FF$ is spanned if and only if it 
is $\text{Spec}(\CC)$-spanned. If $\FF$ is a line bundle, then 
$|\FF|$ separates points if and only if $\FF$ is 
$\text{Spec}(\CC\times\CC)$-spanned. If, in addition, $Z$ is complete, 
$|\FF|$ is very ample if and only if $\FF$ is $\Sigma$-spanned 
for each scheme $\Sigma$ of length 2.

\begin{proposition}\label{verlinde_gg}\setcounter{equation}{0}\  \newline 
For each $k>0$, each $m>kh$, and each scheme $\Sigma$ of length $k$, 
the bundle $V_m$ is $\Sigma$-spanned.
\end{proposition}

\begin{proof}\  \newline 
Let $a$ and $b$ be nonnegative integers such that 
$a+b\leq k$. We claim that
	$$H^i(V_m\otimes\II_{T,J}\otimes\OO_J(-\Theta_{N_1})\otimes\cdots
	\otimes	\OO_J(-\Theta_{N_b}))=0$$
for each $i>0$, each subscheme $T\subseteq J$ of length $a$ and each 
line bundles $N_1,\dots,N_b$ of degree $g-1$. 
We will prove our claim by induction on $a$.

Set $G:=V_m\otimes\OO_J(-\Theta_{N_1})\otimes\cdots\otimes
\OO_J(-\Theta_{N_b})$. Suppose first that $a=0$. 
Since $r_J$ is a finite \'etale covering, ${r_J}_*\OO_J$ is a 
bundle containing $\OO_J$ as a subbundle, and so 
$H^i(F)\subseteq H^i(r_J^*F)$ for each coherent sheaf $F$ on $J$. 
Now, for any line bundle $N$ of degree $g-1$ on $X$, the divisors 
$\Theta_{N_j}$ and $\Theta_N$ are algebraically equivalent. In addition, 
$r_J^*\Theta_N$ is numerically equivalent to $r^2\Theta_N$. So, 
by Formula~(\ref{mult_r}),
	$$r_J^*G\cong\bigoplus_{i=1}^{v_m}\OO_J(D),$$
where $D$ is numerically equivalent to $(mrr_1-br^2)\Theta_N$. Since 
	$$mrr_1>khrr_1=kr^2\geq br^2,$$
the divisor $D$ is ample. Thus $H^i(\OO_J(D))=0$ for all $i>0$, and hence 
$H^i(G)=0$ for all $i>0$.

Suppose now that $a>0$. Let $S\subseteq T$ be a closed subscheme 
of colength 1. By the induction hypothesis,
	\begin{equation}\label{ind-van}
	H^i(G\otimes\II_{S,J})=0\text{ for each $i>0$.}
	\end{equation}
Tensor the natural exact sequence of ideal sheaves 
$0\to\II_{T,J}\to\II_{S,J}\to\II_{S,T}\to 0$ by $G$, and consider the 
resulting long exact sequence in cohomology. 
Since the support of $\II_{S,T}$ is finite, (\ref{ind-van}) yields
	$$H^i(G\otimes\II_{T,J})=0\text{ for each $i>1$.}$$
Furthermore, $H^1(G\otimes\II_{T,J})=0$ if and only if each section of 
$G\otimes\II_{S,T}$ lifts to one of $G\otimes\II_{S,J}$. 
Since $\II_{S,T}$ has length 1, this will be the case if 
$G\otimes\II_{S,J}$ is globally generated. Now, it follows from the 
$M$-regularity 
criterion of \cite{pp} (more precisely from the Corollary on p.~286), 
that $G\otimes\II_{S,J}$ is globally 
generated if $H^i(G\otimes\II_{S,J}\otimes\OO_J(-\Theta_N))=0$ for each 
line bundle $N$ of degree $g-1$ on $X$ and each $i>0$. 
As this is the case by the induction 
hypothesis, we get $H^1(G\otimes\II_{T,J})=0$, thus proving our claim.

Now, let $T\subseteq J$ be a finite subscheme of length at most $k$. 
Tensor the natural exact sequence 
$$0\to\II_{T,J}\to\OO_J\to\OO_T\to 0$$ 
by $V_m$, and consider the resulting exact sequence in cohomology. 
By the claim, $H^1(V_m\otimes\II_{T,J})=0$. So the restriction map 
$H^0(V_m)\to H^0(V_m|_T)$ is surjective.
\end{proof}

The use of Verlinde bundles allows us to reduce the very
ampleness problem on $U_X(r,d)$ to that on $SU_X(r,d)$.

\begin{theorem}\label{very_ampleness}\  \newline
Let $\Sigma$ be a finite scheme of length $k$. If $\LL^m$ is 
$\Sigma$-spanned on $SU_X(r,d)$, and $m>kh$, then 
$\OO_U(m\Theta_F)$ is $\Sigma$-spanned on $U_X(r,d)$.
\end{theorem}

\begin{proof}\  \newline 
We may assume that $\Sigma\subseteq U$. Put 
$S:=\pi_L(\Sigma)$. Note that $S$ has length at most $k$, say $\ell$. 
We need only to prove the following two statements.

\emph{Statement 1:} 
Each section of $\OO_U(m\Theta_F)|_{\Sigma}$ lifts to one of 
$\OO_U(m\Theta_F)|_{\pi_L^{-1}(S)}$. Indeed, since $r_J$ is an \'etale 
covering, $S=r_J(S')$ for some subscheme $S'\subseteq J$ of length $\ell$, 
and thus $\Sigma=\tau(\Sigma')$ for some 
subscheme $\Sigma'\subseteq SU_X(r,L)\times S'$ of length $k$. So we need 
to show that each section of $\tau^*\OO_U(m\Theta_F)|_{\Sigma'}$ 
lifts to one of $\tau^*\OO_U(m\Theta_F)|_{SU_X(r,L)\times S'}$. Now, 
$\tau^*\OO_U(m\Theta_F)=\LL^m\boxtimes\OO_J(mrr_1\Theta_N)$ by 
Formula~(\ref{pullback}). We are left then with proving that the natural map
	$$\lambda\:H^0(\LL^m)\otimes H^0(\OO_{S'})\to 
	H^0(\LL^m\boxtimes {\OO_{S'}}|_{\Sigma'})$$
is surjective. Let $T$ be the image of $\Sigma'$ in $SU_X(r,L)$ 
under the projection. 
Since $\Sigma'$ is isomorphic to $\Sigma$, and $\LL^m$ is $\Sigma$-spanned, 
the natural 
map $H^0(\LL^m)\to H^0(\LL^m|_{T})$ is surjective. Then  
$H^0(\LL^m)\otimes H^0(\OO_{S'})$ surjects onto 
	$$H^0(\LL^m|_{T})\otimes H^0(\OO_{S'})=
	H^0(\LL^m\otimes {\OO_{S'}}|_{T\times S'}).$$
Now, since $\Sigma'\subseteq T\times S'$, each section of 
$\LL^m\otimes {\OO_{S'}}|_{\Sigma'}$ lifts to one of 
$\LL^m\otimes {\OO_{S'}}|_{T\times S'}$. So $\lambda$ is surjective, proving 
the statement.

\emph{Statement 2:} Each section of $\OO_U(m\Theta_F)|_{\pi_L^{-1}(S)}$ 
lifts to one of $\OO_U(m\Theta_F)$. Indeed, 
	$$H^0(\OO_U(m\Theta_F))=H^0(V_m)\text{ and }
	H^0(\OO_U(m\Theta_F)|_{\pi_L^{-1}(S)})=H^0({V_m}|_{S}),$$ 
the latter because the formation of ${\pi_L}_*\OO_U(m\Theta_F)$ commutes 
with base change. 
Since $S$ has length at most $k$, the statement follows from 
Proposition~\ref{verlinde_gg}.
\end{proof}

The proof of Theorem~\ref{very_ampleness} shows that, if $\LL^m$ is 
$\Sigma$-spanned on the stable locus $SU_X^s(r,d)$, 
and if $m>kh$, then $\OO_U(m\Theta_F)$ is $\Sigma$-spanned on $U_X^s(r,d)$.

\begin{proof} (of Theorem B)\  \newline 
Apply Theorem A and Theorem~\ref{very_ampleness}, first with 
$\Sigma:=\text{Spec}(\CC\times\CC)$ and then with 
$\Sigma:=\text{Spec}(\CC[\epsilon]/(\epsilon^2))$, 
keeping in mind the observation after Theorem~\ref{very_ampleness}.
\end{proof}

\begin{proof} (of Theorem C)\  \newline 
If $d$ is odd, then $|\LL|$ is very ample by the main result of \cite{bv2}, 
whereas if $d$ is even and $C$ is nonhyperelliptic, then 
$|\LL|$ is very ample by \cite{vgi}, Thm. 1, p. 134. Now apply 
Theorem~\ref{very_ampleness}, observing that $h=1$ in the first case, 
while $h=2$ in the second case.
\end{proof}

\end{document}